\documentclass[11pt]{article}

\usepackage{amsfonts}

\newtheorem{Theorem}{Theorem}[section]
\newtheorem{Lemma}[Theorem]{Lemma}
\newtheorem{Proposition}[Theorem]{Proposition}
\newtheorem{Corollary}[Theorem]{Corollary}
\newtheorem{Def1}{Definition}[section]
\newtheorem{Exe1}{Example}[section]

\newenvironment{Exe}{\begin{Exe1}\rm}{\end{Exe1}}

\newcommand{\qed}{\hfill $\Box$ \hfill \\}
\long\def\symbolfootnote[#1]#2{\begingroup%
\def\thefootnote{\fnsymbol{footnote}}\footnote[#1]{#2}\endgroup} 

\title{\sc Homomorphisms of complexes via homologies}
\author{Intan Muchtadi-Alamsyah}
\date{}

\begin{document}
\maketitle

\begin{abstract}
\noindent
By Rickard's work, two rings are derived equivalent if there is a tilting complex, constructed from projective modules over the first ring such that the second ring is the endomorphism ring of this tilting complex.

 In this work I describe, under some conditions, the homomorphism space in the derived category between two complexes as iterated pull-back of homomorphism spaces between the homologies of the complexes.
\end{abstract}

\section*{Introduction}

In 1989, J.Rickard \cite{Rickard1} and B.Keller \cite{Keller} have given a necessary and sufficient criterion for the existence of derived equivalences between two rings. Rickard's theorem says that for two rings $\Lambda$ and $\Gamma$ the derived categories  $D^b(\Lambda)$ and $D^b(\Gamma)$ of $\Lambda$ and $\Gamma$ are equivalent as triangulated categories if and only if there exists an object $T$ in $D^b(\Lambda)$, named tilting complex, satisfying similar proprieties as those of a progenerator and such that $\Gamma$ is isomorphic to the endomorphism ring of $T$ in $D^b(\Lambda).$

To determine a priori the endomorphism ring $\Gamma$ of a tilting complex over a given ring $\Lambda$ is difficult. Hence it seems to be important to get an overview of all possible $\Gamma$ by a simpler way than to calculate the endomorphism ring directly.

One possibility is to try to describe the endomorphism ring of tilting complexes via the endomorphism ring of the homologies of the complexes. The aim of this paper is to describe the endomorphism ring of complexes, more generally, the homomorphism space between two complexes as pull-back of homomorphism spaces between their homologies. 

S.K\"onig and A.Zimmermann \cite{Koenig-Zimmermann1} have used this approach to describe the endomorphism ring of a 2-term tilting complex with torsion free homology over a Gorenstein order $\Lambda$ over a complete local Dedekind domain $R$ (see \cite[\S 37]{Curtis-Reiner} for the relevant definitions).  

Klaus Lux, in a private communication with A.Zimmermann, asked whether there exists a generalization of the above result for n-term complexes. We present here the description of $Hom_{D^b(\Lambda)}(T,S)$ for two $n$-term complexes $T$ and $S$ over an algebra which is not necessarily an order and the homologies need not be torsion free. 

For $n = 2,$ we generalize K\"onig and Zimmermann's result. In the absence of their assumpsions, the mapping $\psi^i : Hom_{D^b(\Lambda)}(T,S) \rightarrow Hom_{\Lambda}(H^{i+1}T,H^{i+1}S)$ is not surjective in general; however, the image of $\psi^i$ can be described explicitly. This will give us explicit pull-backs.

The case where $Hom_{D^b(\Lambda)}(T,S)$ is $R$-torsion free is interesting. Even though we don't assume that the right-most homologies are $R$-torsion free, we can replace them by their torsion free parts if the homologies are torsion modules except at the ends of the complexes. In this case $Hom_{D^b(\Lambda)}(T,S)$ is the pull-back of $Im(\psi^0)$ and $Hom_{D^b(\Lambda)}(H^nT/tors, H^nS/tors).$ In particular for $S = T$ a tilting complex and $\Lambda$ a symmetric order, $Hom_{D^b(\Lambda)}(T,T)$ is $R$-torsion free by \cite[Theorem 1]{Zimmermann}. Recall that a symmetric order is an order which is a symmetric algebra \cite{Zimmermann}.

The paper is organized as follows. In section 1 we give the definition of derived categories and we review the Rickard's Morita Theory for derived equivalence. In section 2, we describe $Hom_{D^b(\Lambda)}(T,S)$ as iterated pull-back of homomorphism spaces between their homologies. The image of $\psi^i$ and more explicit pull-backs will be explained in section 3. We will explain the torsion free $Hom$-spaces case in section 4. At the very end we give examples where the main result applies. \\  

{\bf Acknowledgement} I wish to thank Alexander Zimmermann for many discussions and for correcting this paper and also for having introduced me to derived categories. 

\section{Review of derived categories}

\subsection{Category of complexes}

Let $\Lambda$ be an associative ring with 1 and denote by $\Lambda-mod$ the category of finitely generated left $\Lambda$-modules. A {\it complex} $X$ in $\Lambda-mod$ is a sequence of finitely generated $\Lambda$-modules $(X^k)_{k \in \mathbb{Z}}$ and morphisms $(\alpha^k_X)_{k \in \mathbb{Z}}$ in $(Hom_{\Lambda}(X^{k}, X^{k+1}))_{k \in \mathbb{Z}}$
\[ \cdots \rightarrow X^{k-1} \stackrel{\alpha^{k-1}}{\longrightarrow} X^k \stackrel{\alpha^k}{\longrightarrow} X^{k+1} \rightarrow \cdots \]
such that $\alpha^{k-1} \alpha^k = 0$ for all $k \in \mathbb{Z}.$

For two complexes $X$ and $Y,$ a morphism of complexes $\varphi : X \rightarrow Y$ is a sequence $(\varphi^k)_{k \in \mathbb{Z}}$ of morphisms $\varphi^k : X^k \rightarrow Y^k$ such that 
\[ \varphi^k \alpha_Y^k = \alpha_X^k \varphi^{k+1} \] for all $k \in \mathbb{Z}.$ With these definitions we form the {\it category of complexes} of finitely generated $\Lambda$-modules $\mathcal{C}(\Lambda).$

The {\it shift functor} $[i]$ is an endofunctor from $\mathcal{C}(\Lambda)$ to itself where the image of a complex $X$ is another complex $X[i]$ such that
\[ X[i]^k = X^{k+i} \;\; \hbox{and} \;\; \alpha_{X[i]}^k = (-1)^i\alpha_X^{k+i}.\]

For a complex $X$ in $\mathcal{C}(\Lambda)$ we define
\[ H^k(X) = Ker(\alpha^k_X)/Im(\alpha^{k-1}_X) \] and we called the {\it k-th homology} of $X.$ We note that \[ H^k(X) \subseteq X^k/Im(\alpha_X^{k-1}) = Coker(\alpha_X^{k-1}) \;\; \hbox{for all} \; k \in \mathbb{Z}.\]

The category $\Lambda-mod$ embeds into $\mathcal{C}(\Lambda)$ by identifying a finitely generated $\Lambda$-module $M$ with a complex with homology $M$ concentrated in degree 1. We denote by $\mathcal{C}^b(\Lambda)$ the full subcategory of $\mathcal{C}(\Lambda)$ consisting of bounded complexes, i.e the complexes $X$ where $X^k = 0$ for all $k << 0$ and all $k >> 0.$ 

\subsection{Derived categories and derived equivalences}

The objects of $D^b(\Lambda),$ the {\it derived category} of $\Lambda,$ are complexes
\[ P = (\cdots \rightarrow P^k \stackrel{\alpha_P^k}{\longrightarrow} P^{k+1} \rightarrow \cdots ) \]
of finitely generated projective left $\Lambda$-modules $P^k,$ $k \in \mathbb{Z}$ such that 
\[ P^k = 0 \;\hbox{for all} \; k >> 0 \;\; \hbox{and} \;\; H^k(P) = 0 \; \hbox{for all} \; k << 0. \] If $P$ and $Q$ are such complexes, morphisms $\varphi$ and $\psi$ from $P$ to $Q$ are said to be {\it homotopic} if there is some family of morphisms $h^k : P^k \rightarrow Q^{k-1},$ $ k \in \mathbb{Z}$ such that
\[ h^k\alpha_Q^{k-1} + \alpha_P^kh^{k+1} = \varphi^k - \psi^k \;\; \hbox{for all} \;\; k \in \mathbb{Z}. \] This is an equivalence relation, and, by definition the equivalence classes form the morphisms from $P$ to $Q$ in $D^b(\Lambda).$

The category $\Lambda-mod$ (resp. $\mathcal{C}^b(\Lambda))$ embeds into $D^b(\Lambda)$ by choosing for each module (resp. each bounded complex) a fixed projective resolution, which  can then represent the module (resp. the complex) in $D^b(\Lambda).$  

For two complexes $X$ and $Y$ in  $D^b(\Lambda)$ and $k \in \mathbb{Z},$ we identify \\ $Hom_{D^b(\Lambda)}(X,Y[k])$ with $Ext^k_{\Lambda}(X,Y)$ (\cite[Lemma 2]{Keller3}). 

The derived category $D^b(\Lambda)$ is not necessarily an exact category. However, it is triangulated (\cite[Lemma 2.2.3]{Koenig-Zimmermann2}). To each $X$ and $Y$ in $D^b(\Lambda)$ and each morphism $\varphi : X \rightarrow Y$ in $D^b(\Lambda),$ we associate the
 {\it mapping cone} $(C(\varphi), \alpha_{C(\varphi)})$ where for all $k \in \mathbb{Z},$
\begin{eqnarray*} C(\varphi)^k & := &  X^{k+1} \oplus Y^k \\ 
                  \alpha_{C(\varphi)}^k & := & \left( \begin{array}{cc} -\alpha^k_{X[1]} & \varphi^{k+1} \\ 0 & \alpha^k_Y \end{array} \right)\end{eqnarray*}
and the {\it standard triangle}
\[ X \stackrel{\varphi}{\rightarrow} Y \rightarrow C(\varphi) \rightsquigarrow X[1].\]
A {\it triangle} in $D^b(\Lambda)$ is a sequence  \[ X \stackrel{\varphi}{\rightarrow} Y \rightarrow Z \rightsquigarrow X[1] \] isomorphic to a standard triangle, i.e. such that $Z$ is isomorphic to $C(\varphi)$ in $D^b(\Lambda)$ \cite[Proposition 2.5.3]{Koenig-Zimmermann2}. 

Denote by $\Lambda-per$ the full subcategory of $D^b(\Lambda)$ consisting  of the {\it perfect complexes}, i.e  complexes $P$ of finitely generated projective modules $P^k$ where $P^k = 0$ for all $k << 0$ and all $k >> 0.$ Rickard \cite{Rickard1} and Keller \cite{Keller} have given a necessary and sufficient criterion for the existence of a derived equivalence between two rings $\Lambda$ and $\Gamma.$

\begin{Theorem} (Rickard \cite{Rickard1} and Keller \cite{Keller})
The derived categories $D^b(\Lambda)$ and $D^b(\Gamma)$ are equivalent as triangulated categories if and only if there is a complex $T$ in $\Lambda-per$ such that
\begin{enumerate}
\item $Hom_{D^b(\Lambda)}(T,T[i]) = \left\{ \begin{array}{cc}  0 & \hbox{if} \;\; i \neq 0, \\ \Gamma & \hbox{if} \;\; i = 0. \end{array} \right. $
\item the smallest triangulated subcategory generated by direct summands
of finite direct sums of copies of $T,$ inside $\Lambda-per,$ contains $\Lambda$ as a complex concentrated in degree 1.
\end{enumerate}
\end{Theorem}
A complex $T$ satisfying the condition 1 for $i \neq 0$ and the condition 2 is called a {\it tilting complex}  for $\Lambda.$

\section{Description of $Hom_{D^b(\Lambda)}(T,S)$ via homologies}

Let $R$ be a commutative ring and $\Lambda$ an $R$-algebra. Let $T$ and $S$ be complexes \[ \begin{array}{c} T = ( 0 \stackrel{\alpha^0}{\longrightarrow} P^1 \stackrel{\alpha^1}{\longrightarrow} P^2 \stackrel{\alpha^2}{\longrightarrow} \cdots \stackrel{\alpha^{n-1}}{\longrightarrow} P^n \rightarrow 0) \;\; \hbox{and} \;\; \\
S = (0 \stackrel{\beta^0}{\longrightarrow} Q^1 \stackrel{\beta^1}{\longrightarrow} Q^2 \stackrel{\beta^2}{\longrightarrow} \cdots \stackrel{\beta^{n-1}}{\longrightarrow} Q^n \rightarrow 0) \end{array}\]
of projective $\Lambda$-modules $P^i$ and $Q^i$, with homologies concentrated in \\ degree $1,\cdots,n.$ For all $i \in \{ 1,\cdots,n-1 \},$ we define $T^i$ and $S^i$ to be the following quotient-complexes of $T$ and $S$ respectively
\begin{eqnarray*} \lefteqn{ T^i = (  0 \rightarrow Im(\alpha^i) \rightarrow P^{i+1} \stackrel{\alpha^{i+1}}{\longrightarrow} \cdots \stackrel{\alpha^{n-1}}{\longrightarrow} P^n \rightarrow  0)}  \\ & & S^i = ( 0 \rightarrow Im(\beta^i) \rightarrow Q^{i+1} \stackrel{\beta^{i+1}}{\longrightarrow} \cdots \stackrel{\beta^{n-1}}{\longrightarrow} Q^n \rightarrow  0). \end{eqnarray*}
 We denote by $T^0 := T,$ $S^0 := S.$

\begin{Theorem} \label{principal} If for all  $i \in \{ 1, \cdots, n-1 \},$ \begin{eqnarray} \label{condition} \lefteqn{Hom_{\Lambda}(Coker(\alpha^i),Im(\beta^i)) = Hom_{\Lambda}(Im(\alpha^i),Coker(\beta^i)) =} \nonumber \\ & & = Ext^1_{\Lambda}(Coker (\alpha^i),H^iS) = Hom_{\Lambda}(Coker (\alpha^i),Q^i) = 0,\end{eqnarray}  then the diagram 
\[\begin{array}{ccc}  Hom_{D^b(\Lambda)}(T^i,S^i) & \rightarrow &  Hom_{D^b(\Lambda)}(T^{i+1},S^{i+1})          \\ \\
  \phantom{\psi^i} \downarrow \; _{\psi^i}          &             &   \downarrow 
                         \\ \\
Hom_{\Lambda}(H^{i+1}T,H^{i+1}S)      & \rightarrow & Hom_{D^b(\Lambda)}(Coker (\alpha^{i+1}),H^{i+1}S[2])    \end{array} \] is a pull-back diagram for each $i \in \{0, 1,\cdots, n-2\}.$ 
\end{Theorem}
{\bf Remarks} \begin{enumerate} 
\item The condition (\ref{condition}) holds for $i = 0$ because $\alpha^0 = \beta^0 = 0$ and $Coker(\alpha^0) = P^1$ is projective. 
\item Recall that $T^0 = T$ and $S^0 = S,$ so that the theorem presents $Hom_{D^b(\Lambda)}(T,S)$ as an iterated pull-back of homomorphism spaces between $H^{i+1}T$ and $H^{i+1}S$ for all $i \in \{0,1, \cdots, n-1\}.$ 
\item In general, the mapping $\psi^i$ is not surjective. We will determine the image of $\psi^i$ in section 2. 
\item If $R$ is a Dedekind domain with $K = Frac(R)$  and $\Lambda$ is an $R$-order, i.e. an $R$-algebra finitely generated projective as $R$-module such that $K \otimes_R \Lambda$ is a semisimple $K$-algebra, then  \[Hom_{\Lambda}(Coker (\alpha^i),Q^i) = 0\;\; \hbox{ implies} \;\; Hom_{\Lambda}(Coker (\alpha^i),Im(\beta^i)) = 0.\]
\end{enumerate}

We abbreviate for any two complexes $X$ and $Y$ of $\Lambda$-modules \\
$(X,Y) := Hom_{D(\Lambda)}(X,Y),$ $Ext^i(X,Y) = Ext^i_{\Lambda}(X,Y)$ and \[\varphi \circ (X, Y ) := \{\varphi \delta : \delta \in (X, Y )\}. \]
For all $i \in \{0,1, \cdots, n-1\},$ we denote by $L^i := Coker(\alpha^i)$ and by $\iota^i$ the inclusion $H^{i+1}T \hookrightarrow L^i.$ \\

\begin{Lemma}
For all $i \in \{0,1,\cdots, n-2\}$ there exist triangles in $D^b(\Lambda)$
\[ T^i \rightarrow T^{i+1} \rightarrow H^{i+1}T[-i+1] \rightsquigarrow T^i[1] \;\; \hbox{and}\]  \[ S^i \rightarrow S^{i+1} \rightarrow H^{i+1}S[-i+1] \rightsquigarrow S^i[1]\] where $T^0$ is $T$ and $S^0$ is $S.$
\end{Lemma}

\noindent {\bf Proof} We will show that \[ H^{i+1}T[-i] \stackrel{(\hat{\iota})^i}{\longrightarrow} T^i \rightarrow T^{i+1} \rightsquigarrow H^{i+1}T[-i+1] \] is a triangle in $D^b(\Lambda).$  \\
The map $(\hat{\iota})^i : H^{i+1}T[-i] \rightarrow T^i$ is induced by the inclusion $H^{i+1}T \stackrel{\iota^i}{\hookrightarrow} Coker(\alpha^i).$ The mapping cone of $(\hat{\iota})^i$ is the following complex :
\[ 0 \rightarrow Ker (\alpha^{i+1}) \rightarrow P^{i+1} \stackrel{\alpha^{i+1}}{\longrightarrow} \cdots \stackrel{\alpha^{n-1}}{\longrightarrow} P^n \rightarrow 0 \] which is isomorphic to $T^{i+1}.$ Hence, the above sequence is a triangle. \\
Similarly, $H^{i+1}S[-i] \rightarrow S^i \rightarrow S^{i+1} \rightsquigarrow H^{i+1}S[-i+1]$ is a triangle in $D^b(\Lambda).$ The result follows from the axiom (TR3) in \cite[Theorem 2.3.1]{Koenig-Zimmermann2}. \qed

From now on we assume that for all $i \in \{ 1, \cdots, n-1 \},$ \begin{eqnarray} \label{condition2} Hom_{\Lambda}(Im(\alpha^i),Coker(\beta^i)) &  =  Hom_{\Lambda}(L^i,Im(\beta^i))  = 0 \nonumber \\  \hbox{and} \;\;  Hom_{\Lambda}(L^i,Q^i)  = 0. \end{eqnarray} The assumption $Ext^1(L^i,H^iS) = 0$ will be only needed for Lemma \ref{thetai}. \\

Let us fix an $i \in \{0, \cdots, n-2 \}$ for the rest of this section and if $i > 0$ we denote by
\[\cdots \rightarrow \bar{P}^0 \rightarrow \bar{P}^1 \rightarrow \cdots \rightarrow \bar{P}^{i-1} \rightarrow 0\]
the first terms of a projective resolution of $Ker(\alpha^i)[-i+2],$  and by 
\[\cdots \rightarrow \tilde{P}^0 \rightarrow \tilde{P}^1 \rightarrow \cdots \rightarrow \tilde{P}^{i-1} \rightarrow \tilde{P}^i  \rightarrow 0\]
the first terms of a projective resolution of $Ker(\alpha^{i+1})[-i+1].$ \\

We use similar notation for projective resolutions of $Ker(\beta^i)[-i+2]$ and $Ker(\beta^{i+1})[-i+1]$ (with $P$ replaced by $Q$). 

\begin{Lemma} \label{Ti}  \begin{enumerate} \item We have $(T^i,H^{i+1}S[-i+1]) = 0$ and $(T^i,S^{i+1}[-1]) = 0.$ 
 \item There exists a short exact sequence \[0 \rightarrow (T^i,H^{i+1}S[-i]) \rightarrow (T^i, S^i) \rightarrow (T^i, S^{i+1}) \rightarrow 0.\] \end{enumerate}
 \end{Lemma} 
 
\noindent {\bf Proof} We have a triangle \[S^i \rightarrow S^{i+1} \rightarrow H^{i+1}S[-i+1] \rightsquigarrow S^i[1]\] to which we apply $(T^i,-)$ and 
which then gives rise to a long exact sequence part of which looks as follows: 
 \begin{eqnarray} \label{les1} \lefteqn{ \cdots \rightarrow (T^i, S^{i+1}[-1]) \rightarrow (T^i, H^{i+1}S[-i]) \rightarrow} \nonumber\\ & &  \rightarrow (T^i,S^i) \rightarrow (T^i, S^{i+1}) \rightarrow (T^i, H^{i+1}S[-i+1]) \rightarrow \cdots \end{eqnarray} 
\noindent {\it Proof of 1 :} \begin{enumerate} \item[(a)] $(T^i, H^{i+1}S[-i+1]) = 0$ \\
For $i = 0,$ the result follows from the fact that  $T^0$ is concentrated in degrees $1,\cdots,n,$ whereas the projective resolution of $H^1S[1]$ is $0$ in degrees $1,\cdots,n.$\\
For $i > 0,$ a morphism in $(T^i, H^{i+1}S[-i+1])$ is given by a commutative diagram 
\[ \begin{array}{ccccccccccccc} 
\cdots & \rightarrow & \bar{P}^{i-1} & \rightarrow & P^i        & \rightarrow & P^{i+1}    & \rightarrow & \cdots & \rightarrow & P^n & \rightarrow & 0 \\
       &             & \downarrow    &             & \downarrow &             & \downarrow &             &        &             &     &             &   \\
       &             & 0             & \rightarrow & H^{i+1}S   & \rightarrow & 0          &             &        &             &     &             &
\end{array} \]The morphism $P^i \rightarrow H^{i+1}S$ factors through the cokernel of \\ $(\bar{P}^{i-1} \rightarrow P^i)$ which is $Im(\alpha^i).$  Hence we get \[(T^i, H^{i+1}S[-i+1]) = \frac{(Im(\alpha^i),H^{i+1}S)}{j^i \circ (P^{i+1},H^{i+1}S)}\] where $j^i$ is the inclusion $Im(\alpha^i) \hookrightarrow P^{i+1}$ and \[j^i \circ (P^{i+1},H^{i+1}S) = \{j^i \varphi : \varphi \in (P^{i+1},H^{i+1}S)\} \subseteq (Im(\alpha^i),H^{i+1}S).\] Since $(Im(\alpha^i),H^{i+1}S) \subseteq (Im(\alpha^i),Coker(\beta^i)) = 0$ we get \\ $(T^i, H^{i+1}S[-i+1]) = 0.$ 

\item[(b)] $(T^i,S^{i+1}[-1]) = 0.$ 
  
 Any such morphism is given by a commutative diagram as follows :
 {\footnotesize
  \[\begin{array}{ccccccccccccc}
\cdots & \rightarrow  & P^{i+1}    & \stackrel{\alpha^{i+1}}{\rightarrow} & P^{i+2}          & \stackrel{\alpha^{i+2}}{\rightarrow} &  \cdots & \rightarrow & P^n        & \rightarrow                          & 0          &             & \\
       &              & \downarrow &                                          & \downarrow       &                                     &   & & \downarrow &                                                                                         & \downarrow &             & \\
       &              & 0          & \rightarrow                          & Im(\beta^{i+1}) & \rightarrow  & \cdots & \rightarrow & Q^{n-1}    & \stackrel{\beta^{n-1}}{\rightarrow} & Q^n        & \rightarrow & 0  \end{array}\] }
 
As the left-most square commutes, the morphism $P^{i+2} \rightarrow Im(\beta^{i+1})$ factors through $L^{i+1}$. By assumption (\ref{condition2}), $(L^{i+1},Im(\beta^{i+1})) = 0,$ hence we may assume that the homomorphism $P^{i+2} \rightarrow Im(\beta^{i+1})$ is $0$. \\
Now since $(L^j,Q^j) = 0$ for all $j \in \{ i+2,\cdots, n-1\},$ we can apply an analogous argument to the homomorphisms  $P^{j+1} \rightarrow Q^j$ for all $j \in \{i+2,\cdots, n-1\}.$ We therefore obtain $(T^i, S^{i+1}[-1]) = 0.$ 
 \end{enumerate}
\noindent{\it Proof of 2} : The result follows immediately from 1 applied to the long exact sequence (\ref{les1}). \qed

We denote by $\cdots \rightarrow \hat{P}^{i-2} \rightarrow \hat{P}^{i-1} \rightarrow \tilde{P}^i \rightarrow 0$  the first terms of a projective resolution of $H^{i+1}T[-i+1].$ 

\begin{Lemma} \label{Hi} We have \begin{enumerate} \item $(H^{i+1}T[-i+1],S^{i+1}) = 0,$
 \item $(H^{i+1}T[-i],S^{i+1}) = 0,$ \item $(T^i,S^{i+1}) \cong (T^{i+1},S^{i+1}).$ \end{enumerate}
\end{Lemma}

\noindent {\bf Proof}  We apply  $(-,S^{i+1})$ to the triangle $$T^i \rightarrow T^{i+1} \rightarrow H^{i+1}T[-i+1] \rightsquigarrow T^i[1]$$ and
  get a long exact sequence part of which looks as follows:\begin{eqnarray}  \label{les2} \lefteqn{\cdots \rightarrow (H^{i+1}T[-i+1],S^{i+1}) \rightarrow (T^{i+1},S^{i+1}) \rightarrow} \nonumber\\ & & \rightarrow (T^i,S^{i+1}) \rightarrow (H^{i+1}T[-i],S^{i+1}) \rightarrow \cdots \end{eqnarray}
 
 \begin{enumerate}
 \item $(H^{i+1}T[-i+1],S^{i+1}) = 0.$

Given such a morphism, we get a commutative diagram
\[\begin{array}{ccccccccccccc}
\cdots & \rightarrow & \hat{P}^{i-1}     & \rightarrow & \tilde{P}^i        & \rightarrow & 0          &             &        &             &     &             &     \\
       &             & \downarrow &             & \downarrow &             & \downarrow &             &        &             &     &             &  \\
\cdots & \rightarrow & \tilde{Q}^{i-1}     & \rightarrow & \tilde{Q}^i        & \rightarrow & Q^{i+1}        & \rightarrow & \cdots & \rightarrow & Q^n & \rightarrow & 0 \end{array} \]
 Now, the morphism $\tilde{P}^i \rightarrow \tilde{Q}^i$ factors through the kernel of the morphism  $\tilde{Q}^i \rightarrow Q^{i+1},$ hence through its projective cover $\tilde{Q}^{i-1}.$ An analogous argument shows that for all $j < i,$ the homomorphism $\hat{P}^j \rightarrow \tilde{Q}^j$ factors through $\tilde{Q}^{j-1}.$ Hence the chain map represented in the diagram above is homotopic to zero.
 
\item $(H^{i+1}T[-i],S^{i+1}) = 0.$

Again, such a morphism is given by a commutative diagram
\[\begin{array}{ccccccccccccc}
\cdots & \rightarrow & \hat{P}^{i-1}     & \rightarrow & \tilde{P}^i        & \rightarrow & 0          &             &        &             &     &             &     \\
       &             & \downarrow &             & \downarrow &             & \downarrow &             &        &             &     &             &  \\
\cdots & \rightarrow & \tilde{Q}^i     & \rightarrow & Q^{i+1}       & \rightarrow & Q^{i+2}        & \rightarrow & \cdots & \rightarrow & Q^n & \rightarrow & 0 \end{array} \]
With the same argument as in 1, the homomorphism $\tilde{P}^i \rightarrow Q^{i+1}$ factors through $\tilde{Q}^i$ and the homomorphism $\hat{P}^j \rightarrow \tilde{Q}^{j+1}$ for all $j < i$ factors through $\tilde{Q^j}.$ Hence the chain map in the diagram above is homotopic to zero.
\item By applying 1 and 2 to the long exact sequence (\ref{les2}), we obtain \\ $(T^i,S^{i+1}) \cong (T^{i+1},S^{i+1}).$ \qed \end{enumerate}
  
\begin{Corollary} \label{pii} From Lemma \ref{Ti} and Lemma \ref{Hi} we get that  \[0 \rightarrow (T^i,H^{i+1}S[-i]) \rightarrow (T^i,S^i) \rightarrow (T^{i+1},S^{i+1}) \rightarrow 0\] is a short exact sequence. 
\end{Corollary}

\begin{Lemma} \label{zetai} \begin{enumerate} \item  $Ext^1(Im(\alpha^{i+1}),H^{i+1}S) = (L^{i+1},H^{i+1}S[2]).$
\item  The mapping $H^{i+1}T \hookrightarrow L^i$ induces a short exact sequence \[0 \rightarrow \iota^i \circ (L^i,H^{i+1}S) \rightarrow (H^{i+1}T,H^{i+1}S) \stackrel{\eta^i}{\rightarrow} (L^{i+1},H^{i+1}S[2]) \rightarrow 0\] \end{enumerate}
\end{Lemma}

\noindent {\bf Proof} Apply $(-,H^{i+1}S)$ to the exact sequence $$0 \rightarrow H^{i+1}T \rightarrow L^i \rightarrow Im(\alpha^{i+1}) \rightarrow 0$$ to get a long exact sequence part of which is :
\begin{eqnarray} \label{les4} \lefteqn{0 \rightarrow (Im(\alpha^{i+1}),H^{i+1}S) \rightarrow (L^i, H^{i+1}S) \rightarrow}  \\ & & \rightarrow (H^{i+1}T,H^{i+1}S) \stackrel{\eta^i}{\rightarrow} Ext^1(Im(\alpha^{i+1}),H^{i+1}S) \rightarrow Ext^1(L^i, H^{i+1}S) \rightarrow \cdots \nonumber\end{eqnarray} 
\noindent {\it Proof of 1} : $Ext^1(Im(\alpha^{i+1}),H^{i+1}S) = (L^{i+1},H^{i+1}S[2])$ since $Im(\alpha^{i+1})$ is the first syzygy of $L^i.$ 

\noindent{\it Proof of 2} : 
\begin{enumerate}
\item[(a)]$Ext^1(L^i, H^{i+1}S) = 0.$ \\
For $i = 0,$ $Ext^1(L^0, H^1S) = 0$ since $L^0 = P^1$ is projective.\\ For $i > 0,$ the hypothesis $(Im(\alpha^i),Coker(\beta^i)) = 0$ implies \\ $(Im(\alpha^i),H^{i+1}S) = 0$ and therefore  $Ext^1(L^i, H^{i+1}S) = 0.$ 
\item[(b)]The kernel of $(H^{i+1}T,H^{i+1}S) \stackrel{\eta^i}{\rightarrow} (L^{i+1},H^{i+1}S[2])$ is the image of  \\ $(L^i, H^{i+1}S) \rightarrow (H^{i+1}T,H^{i+1}S)$ and this one is $\iota^i \circ(L^i,H^{i+1}S).$ 
\end{enumerate}
From 1, (a), (b) and the long exact sequence (\ref{les4}), we obtain the exact sequence $0 \rightarrow \iota^i \circ (L^i,H^{i+1}S) \rightarrow (H^{i+1}T,H^{i+1}S) \stackrel{\eta^i}{\rightarrow} (L^{i+1},H^{i+1}S[2]) \rightarrow 0.$ \qed

\begin{Lemma} \label{Ext1i} \begin{enumerate} \item We have \begin{enumerate} \item[(a)] $(T^{i+1},H^{i+1}S[-i]) = Ext^1(L^{i+1},H^{i+1}S)$ and \item[(b)] $(T^{i+1},H^{i+1}S[-i+1]) = (L^{i+1},H^{i+1}S[2]).$ \end{enumerate}
\item There exists a short exact sequence
\[ 0 \rightarrow Ext^1(L^{i+1},H^{i+1}S) \rightarrow (T^i,H^{i+1}S[-i]) \rightarrow \iota^i \circ(L^i,H^{i+1}S) \rightarrow 0 \] \end{enumerate}
\end{Lemma}

\noindent {\bf Proof} We apply  $(-,H^{i+1}S[-i])$ to the triangle  $$T^i \rightarrow T^{i+1} \rightarrow H^{i+1}T[-i+1] \rightsquigarrow T^i[1]$$ to get a long exact sequence
  \begin{eqnarray} \label{les5} \lefteqn{ \cdots \rightarrow (H^{i+1}T[-i+1],H^{i+1}S[-i]) \rightarrow (T^{i+1},H^{i+1}S[-i]) \rightarrow (T^i,H^{i+1}S[-i]) \rightarrow} \nonumber \\ & & \rightarrow(H^{i+1}T[-i],H^{i+1}S[-i])  
\rightarrow (T^{i+1}[-1],H^{i+1}S[-i]) \rightarrow  \cdots \end{eqnarray} 
\noindent {\it Proof of 1} :
\begin{enumerate}
\item[(a)] $(T^{i+1},H^{i+1}S[-i]) = Ext^1(L^{i+1},H^{i+1}S)$ 

A morphism $T^{i+1} \rightarrow H^{i+1}S[-i]$ gives rise to a commutative diagram
{\footnotesize \[\begin{array}{ccccccccccccccc}
 \cdots & \rightarrow & \tilde{P}^{i-1} & \rightarrow & \tilde{P}^i        & \stackrel{\tilde{\alpha}^i}{\rightarrow} & P^{i+1}    & \stackrel{\alpha^{i+1}}{\rightarrow} & P^{i+2}    & \rightarrow & \cdots & \rightarrow & P^n & \rightarrow & 0 \\
        &             &               &             & \downarrow &             & \downarrow &             & \downarrow &             &        &             &     &             &   \\
        &             &               &             & 0          & \rightarrow & H^{i+1}S  & \rightarrow & 0          &             &        &             &     &             &   \end{array}\]}
and the morphism $P^{i+1} \rightarrow H^{i+1}S$ factors through $Coker(\tilde{\alpha}^i)$ which is $Im(\alpha^{i+1}).$ Hence, denoting the embedding $Im(\alpha^{i+1}) \hookrightarrow P^{i+2}$ by $j^{i+1},$ we have
\[(T^{i+1},H^{i+1}S[-i]) = \frac{(Im(\alpha^{i+1}),H^{i+1}S)}{j^{i+1} \circ (P^{i+2},H^{i+1}S)} = Ext^1(L^{i+1},H^{i+1}S). \] 

\item[(b)] $(T^{i+1}[-1],H^{i+1}S[-i]) = (T^{i+1},H^{i+1}S[-i+1]) = (L^{i+1},H^{i+1}S[2]).$ 

A mapping $T^{i+1}[-1] \rightarrow H^{i+1}S[-i]$ is given by a commutative diagram
{\footnotesize \[\begin{array}{ccccccccccccccc}
 \cdots & \rightarrow & \tilde{P}^{i-1} & \stackrel{\tilde{\alpha}^{i-1}}{\rightarrow} & \tilde{P}^i        & \stackrel{\tilde{\alpha}^i}{\rightarrow} & P^{i+1}    & \rightarrow & P^{i+2}    & \rightarrow & \cdots & \rightarrow & P^n & \rightarrow & 0 \\
        &             & \downarrow      &             & \downarrow \gamma       &                                         & \downarrow &             &            &             &        &             &     &             &  \\
        &             &     0           & \rightarrow & H^{i+1}          & \rightarrow                             & 0          &             &             &             &       &             &    &              &  \end{array}\]}

Again, $\gamma$ factors through $Coker(\tilde{\alpha}^{i-1}) = Im(\tilde{\alpha}^i).$ Hence, denoting the embedding $Im(\tilde{\alpha}^i) \hookrightarrow P^{i+1} $ by $\rho^i,$ we have : \[ (T^{i+1},H^{i+1}S[-i+1]) = \frac{(Im(\tilde{\alpha}^i),H^{i+1}S)}{\rho^i \circ (P^{i+1},H^{i+1}S)} = Ext^1(Im(\alpha^{i+1}),H^{i+1}S). \] 
From Lemma \ref{zetai}, we get $Ext^1(Im(\alpha^{i+1}),H^{i+1}S) = (L^{i+1},H^{i+1}S[2]).$ \end{enumerate}

\noindent {\it Proof of 2} :
It is clear that \begin{equation} \label{clear} (H^{i+1}T[-i+1],H^{i+1}S[-i]) = Ext^{-1}(H^{i+1}T,H^{i+1}S) = 0. \end{equation}
The image of $(T^i,H^{i+1}S[-i]) \rightarrow (H^{i+1}T,H^{i+1}S)$ is the kernel of  \\ $(H^{i+1}T,H^{i+1}S) \stackrel{\eta^i}{\rightarrow} (L^i,H^{i+1}S[2])$ which is $\iota^i \circ (L^i,H^{i+1}S)$ by Lemma \ref{zetai}. \\
The result follows from 1, (\ref{clear}), and the exactness of (\ref{les5}). \qed
 
 \begin{Lemma}\label{psii}
The morphism $T^i \rightarrow T^{i+1}$ gives rise to an exact sequence \[0 \rightarrow (T^{i+1},S^i) \rightarrow (T^i,S^i) \stackrel{\psi^i}{\rightarrow} (H^{i+1}T,H^{i+1}S).\]
\end{Lemma}

\noindent {\bf Proof} We apply $(-,S^i)$ to the triangle $$T^i \rightarrow T^{i+1} \rightarrow H^{i+1}T[-i+1] \rightsquigarrow T^i[1]$$ to get a long exact sequence part of which is 
\begin{eqnarray*} \lefteqn{ \cdots \rightarrow (H^{i+1}T[-i+1],S^i) \rightarrow (T^{i+1},S^i) \rightarrow (T^i,S^i) \rightarrow} \\ & & \rightarrow (H^{i+1}T[-i],S^i) \rightarrow (T^{i+1}[-1],S^i) \rightarrow \cdots \end{eqnarray*} 
\begin{enumerate}
\item $(H^{i+1}T[-i+1],S^i) = 0.$ \\
For $i = 0,$ the result follows from the fact that  $S^0$ is concentrated in degrees $1,\cdots,n$ whereas the projective resolution of $H^1T[1]$ is $0$ in degrees $1, \cdots,n.$ \\
For $i > 0,$ any morphism in $(H^{i+1}T[-i+1],S^i)$ is given by the commutative diagram as follows : 
\[\begin{array}{ccccccccccccc}
\cdots & \rightarrow & \tilde{P}^{i-1}     & \rightarrow & \tilde{P}^i        & \rightarrow & 0          &             &        &             &     &             &     \\
       &             & \downarrow &             & \downarrow &             & \downarrow &             &        &             &     &             &  \\
\cdots & \rightarrow & \bar{Q}^{i-1}     & \rightarrow & Q^i        & \rightarrow & Q^{i+1}        & \rightarrow & \cdots & \rightarrow & Q^n & \rightarrow & 0 \end{array} \]
With the same argument as in the first part of the proof of Lemma \ref{Hi} applied to $\tilde{P}^i \rightarrow Q^i$ and $\tilde{P}^j \rightarrow \bar{Q}^j,$  for all $j < i,$ we get \[(H^{i+1}T[-i+1],S^i) = 0.\] 

\item $(H^{i+1}T[-i],S^i) = (H^{i+1}T,H^{i+1}S).$ \\
We have a long exact sequence 
\begin{eqnarray*} \lefteqn{ \cdots \rightarrow (H^{i+1}T[-i],S^{i+1}[-1]) \rightarrow (H^{i+1}T[-i],H^{i+1}S[-i]) \rightarrow } \\ & & \rightarrow (H^{i+1}T[-i],S^i) \rightarrow (H^{i+1}T[-i],S^{i+1}) \rightarrow \cdots \end{eqnarray*} coming from applying $(H^{i+1}T[-i],-)$ to the triangle \[S^i \rightarrow S^{i+1} \rightarrow H^{i+1}S[-i+1] \rightsquigarrow S^i[1].\] 
By Lemma \ref{Hi},  $(H^{i+1}T[-i],S^{i+1}[-1]) = 0$ and \\ $(H^{i+1}T[-i],S^{i+1}) = 0.$ Hence, $(H^{i+1}T[-i],S^i) = (H^{i+1}T,H^{i+1}S).$ 
\end{enumerate} 
We obtain an exact sequence  \[0 \rightarrow (T^{i+1},S^i) \rightarrow (T^i,S^i) \stackrel{\psi^i}{\rightarrow} (H^{i+1}T,H^{i+1}S).\] \qed

\noindent {\bf Remark} We will give an example in section 5 (Example \ref{exe4}) which shows that  $(T^i,S^i) \rightarrow (H^{i+1}T,H^{i+1}S)$ is not surjective in general, i.e. the mapping  $(H^{i+1}T[-i],S^i) \rightarrow (T^{i+1}[-1],S^i)$ is not always $0.$ 

\begin{Lemma} \label{thetai}
If $Ext^1(L^{i+1},H^{i+1}S) = 0,$ then the morphism $S^i \rightarrow S^{i+1}$ gives rise to an exact sequence \[0 \rightarrow (T^{i+1},S^i) \rightarrow (T^{i+1},S^{i+1}) \rightarrow (L^{i+1},H^{i+1}S[2]).\]
\end{Lemma}

\noindent {\bf Proof} Applying $(T^{i+1},-)$ to the triangle  $$S^i \rightarrow S^{i+1} \rightarrow H^{i+1}S[-i+1] \rightsquigarrow S^i[1]$$ we get a long exact sequence
\begin{eqnarray*} \lefteqn{ \cdots \rightarrow (T^{i+1},H^{i+1}S[-i]) \rightarrow (T^{i+1},S^i) \rightarrow (T^{i+1},S^{i+1}) \rightarrow} \\ & & \rightarrow (T^{i+1},H^{i+1}S[-i+1]) \rightarrow (T^{i+1},S^i[1]) \rightarrow \cdots \end{eqnarray*} 
From Lemma \ref{Ext1i} and our hypothesis we get \[(T^{i+1},H^{i+1}S[-i]) = Ext^1(L^{i+1},H^{i+1}S) = 0.\] Now
$(T^{i+1},H^{i+1}S[-i+1]) = (L^{i+1},H^{i+1}S[2])$ by Lemma \ref{Ext1i} and the result follows.\qed

\noindent {\bf PROOF OF THEOREM \ref{principal}} :

As a consequence of Corollary \ref{pii}, Lemma \ref{zetai}, Lemma \ref{Ext1i}, Lemma \ref{psii} and Lemma \ref{thetai} we get the following diagram 

{\small \[\begin{array}{ccccccccc}
   &             &  0                        &             &     0                   &             &    0
                             &             &     \\
   &             & \downarrow                &             &   \downarrow            &             &  \downarrow                                      &             &     \\
   &             & Ext^1(L^{i+1},H^{i+1}S)  &             & (T^{i+1},S^i) &           & (T^{i+1},S^i)      &             &    \\
   &             & \downarrow               &             &   \downarrow            &             &  \downarrow                                      &             &     \\
 0 & \rightarrow & (T^i,H^{i+1}S[-i])& \rightarrow & (T^i,S^i) & \rightarrow &  (T^{i+1},S^{i+1})                              & \rightarrow & 0  \\
   &             &         \downarrow       &             &   \phantom{\psi^i} \downarrow \;_{\psi^i}         &             &  \downarrow 
                             &             &    \\
 0 & \rightarrow & \iota^i \circ (L^i,H^{i+1}S)& \rightarrow & (H^{i+1}T,H^{i+1}S)      & \rightarrow & (L^{i+1},H^{i+1}S[2])   & \rightarrow & 0 \\
   &             & \downarrow              &             &                          &             &                         &             &  \\
   &             &    0                    &             &                          &             &                         &             & 
   
 \end{array} \]}

Using the identification $(T^i,S^{i+1}) \cong (T^{i+1},S^{i+1})$ shown in Lemma \ref{Ti} and the hypothesis $(Im(\alpha^i),Coker(\beta^i)) = 0,$ it follows that the bottom right square of the diagram is commutative. Similarly, the bottom left square is commutative.

If  $Ext^1(L^{i+1},H^{i+1}S) = 0,$ then  $(T^i,H^{i+1}S[-i]) \stackrel{\cong}{\rightarrow} \iota^i \circ (L^i,H^{i+1}S)$ and the snake lemma gives the isomorphism $(T^{i+1},S^i) \stackrel{\cong}{\rightarrow} (T^{i+1},S^i)$ in the above diagram.  

This fact and the exactness of the horizontal sequences give us,  for all $i \in \{0,1,\cdots, n-2\},$ the pull-back diagram
 \[\begin{array}{ccc}
 Hom_{D^b(\Lambda)}(T^i,S^i)      & \rightarrow &  Hom_{D^b(\Lambda)}(T^{i+1},S^{i+1})        \\
          \downarrow              &             &   \downarrow            \\
 Hom_{\Lambda}(H^{i+1}T,H^{i+1}S) & \rightarrow & Hom_{D^b(\Lambda)}(L^{i+1},H^{i+1}S[2])    
 \end{array} \] Theorem \ref{principal} is now proved. \qed
 
\noindent {\bf Remark} In fact, the snake lemma shows that $(T^{i+1},S^i) \rightarrow (T^{i+1},S^i)$ is surjective. If $R$ is a noetherian ring, a surjective endomorphism of a noetherian module is an isomorphism. Hence $Ext^1(L^{i+1},H^{i+1}S) = 0.$ But we have already used $Ext^1(L^{i+1},H^{i+1}S) = 0$ to establish the right vertical exact sequence. 

\section{More explicit pull-backs}

In the previous section, we got the pull-back diagrams
\[\begin{array}{ccc}
 Hom_{D^b(\Lambda)}(T^i,S^i)      & \rightarrow &  Hom_{D^b(\Lambda)}(T^{i+1},S^{i+1})        \\ \\
          \phantom{\psi^i} \downarrow \; _{ \psi^i}              &             &  \downarrow            \\ \\
 Hom_{\Lambda}(H^{i+1}T,H^{i+1}S) & \rightarrow & Hom_{D^b(\Lambda)}(L^{i+1},H^{i+1}S[2])    
 \end{array} \]  for all $i \in \{0, 1,\cdots, n-2\}.$ \\
 In this section, we will determine the image of $\psi^i,$ for all $i \in \{0,1,\cdots,n-2\}.$ \\
  
We denote for all $i \in \{0,1,\cdots, n-2\},$ \[A^i := Im(\psi^i) = Ker [(H^{i+1}T,H^{i+1}S) \rightarrow (T^{i+1},S^i[1]) ]. \]  

 For all $i \in \{0,1,\cdots,n-2\},$ we denote by $\iota^i_T$ the inclusion $H^{i+1}T \hookrightarrow Coker(\alpha^i),$  by $\sigma^i_T$ the composition $L^i \rightarrow Im(\alpha^{i+1}) \rightarrow P^{i+2}$ and by $\iota^i_S$ and $\sigma^i_S$ the corresponding mappings for $S$. We remark that $\sigma^0_T = \alpha^1,$ $\sigma^0_S = \beta^1,$ and $\iota^i_T = \iota^i.$ 

\begin{Lemma} For all $i \in \{0,1,\cdots,n-2\},$ we have : \\
$A^i$ is the set of all $\varphi^i \in (H^{i+1}T,H^{i+1}S)$ such that there exists a sequence \[(\gamma^{i+1},\gamma^{i+2},\cdots,\gamma^n) \in (Coker(\alpha^i),Coker(\beta^i)) \times (P^{i+2},Q^{i+2}) \times \cdots \times (P^n,Q^n)\] satisfying 
\begin{enumerate}
\item  $\varphi^i \iota^i_S = \iota^i_T \gamma^{i+1}$, 
\item $\gamma^{i+1} \sigma^i_S = \sigma^i_T \gamma^{i+2},$ and
\item for all $j \in \{i+2,\cdots,n-1\},$  $\gamma^j \beta^j = \alpha^j \gamma^{j+1}.$
\end{enumerate}
\end{Lemma}

\noindent {\bf Proof} 
We denote by $\phi^i$ the projective cover mapping $\tilde{P^i} \rightarrow Ker(\alpha^{i+1}),$ \\ by $\epsilon^i$ the mapping $Ker(\alpha^{i+1}) \rightarrow H^{i+1}T,$ and by $\tilde{\alpha}^i$ the mapping $\tilde{P^i} \rightarrow P^{i+1}$ (we remark that $\epsilon^0 =id$). \\

We have the following diagram :
{\small \[\begin{array}{ccccccccc}
 \cdots & \rightarrow & \tilde{P}^{i-1} & \rightarrow & \tilde{P^i}                        & \stackrel{\tilde{\alpha}^i}{\rightarrow} & P^{i+1}     \stackrel{\alpha^{i+1}}{\rightarrow}  P^{i+2}  \rightarrow   \cdots    \rightarrow & P^{n-1} \rightarrow  P^n   \rightarrow  0 &  T^{i+1}  \\
        &             & \downarrow      &             &  \downarrow \;_{\phi^i \epsilon^i} &                                          & \downarrow \phantom{  \cdots    \rightarrow  P^{n-1} \rightarrow  P^n}  & & \downarrow  \\
        &             & 0               & \rightarrow & H^{i+1}T                           & \rightarrow                              & 0 \phantom{ \cdots    \rightarrow  P^{n-1} \rightarrow  P^n}   &  &  H^{i+1}T  \\
         &             & \downarrow      &             &  \downarrow \;_{\varphi^i} &                                          & \downarrow \phantom{ \cdots    \rightarrow  P^{n-1} \rightarrow  P^n} &  &\downarrow  \\
        &             & 0               & \rightarrow & H^{i+1}S                           & \rightarrow                              & 0 \phantom{ \cdots    \rightarrow  P^{n-1} \rightarrow  P^n}  &  &    H^{i+1}S  \\
         &             & \downarrow      &             &  \downarrow \;_{\iota^i_S} &                                          & \downarrow \phantom{  \cdots    \rightarrow  P^{n-1} \rightarrow  P^n} &    &             \downarrow  \\
  & & 0 & \rightarrow & Coker(\beta^i) & \rightarrow & Q^{i+2} \stackrel{\beta^{i+2}}{\rightarrow}   Q^{i+3}  \rightarrow   \cdots  \rightarrow & Q^n   \rightarrow 0 \phantom{ P^n   \rightarrow  0}&  S^i[1]  \\
 \end{array}\] }

If  $\varphi^i \in A^i$ then $\varphi^i \in (H^{i+1},H^{i+1}S)$ and 
\begin{enumerate}
\item  there exists a mapping $\bar{\gamma}^{i+1} : P^{i+1} \rightarrow Coker(\beta^i)$ such that \[\phi^i \epsilon^i \varphi^i \iota^i_S = \tilde{\alpha}^i \bar{\gamma}^{i+1}, \;\;\hbox{and} \]  
\item there exists a sequence $(\gamma^j)_{j \in \{i+2,\cdots,n-1\}} \;\; \gamma^j : P^j \rightarrow Q^j$ where 
\begin{enumerate} \item $\bar{\gamma}^{i+1} \sigma^i_S = \alpha^{i+1} \gamma^{i+2}$ and \item for all $j \in \{i+2,\cdots,n-1\}$ we have $ \gamma^j \beta^j = \alpha^j \gamma^{j+1}.$ \end{enumerate} \end{enumerate}

Since $(Im(\alpha^i),Coker(\beta^i)) = 0,$ there exists $\gamma^{i+1} \in (Coker(\alpha^i),Coker(\beta^i))$ such that 
\[\bar{\gamma}^{i+1} = \pi^i \gamma^{i+1}, \;\;\varphi^i \iota^i_S = \iota^i_T \gamma^{i+1} \;\; \hbox{and} \;\;\gamma^{i+1} \sigma^i_S = \sigma^i_T \gamma^{i+2},\] where $\pi^i$ is the  mapping $P^{i+1} \rightarrow L^i.$ This finishes the proof of the lemma. \qed

 Hence an element of $A^i$ is a mapping $H^{i+1}T \rightarrow H^{i+1}S$ which induces a homomorphism between the complex \[H^{i+1}T \hookrightarrow Coker(\alpha^i) \rightarrow P^{i+2} \rightarrow \cdots \rightarrow P^n\] and the complex \[H^{i+1}S \hookrightarrow Coker(\beta^i) \rightarrow Q^{i+2} \rightarrow \cdots \rightarrow Q^n.\]
                        
Moreover, for $S = T,$ $A^i$ has a multiplicative stucture as a subring of \\ $(H^{i+1}T,H^{i+1}S)$ for all $i \in \{0,\cdots,n-2\}.$ \\
                        
We denote for all  $i \in \{0,\cdots,n-2\},$  \[\Omega^i := Ker((L^{i+1},H^{i+1}S[2]) \rightarrow (T^{i+1},S^i[1])).\]We remark that if $R$ is a Dedekind domain and $\Lambda$ is an $R$-order, $\Omega^i$ is an $R$-torsion $\Lambda$-module and  $K \otimes_R \Omega^i = 0$ where $K = Frac(R).$

\begin{Lemma} \label{rhoi}
For all $i \in \{0,1,\cdots,n-2\},$ there exists a short exact sequence \\
\[0 \rightarrow \iota^i \circ(L^i,H^{i+1}S) \rightarrow A^i \stackrel{\rho^i}{\rightarrow} \Omega^i \rightarrow 0.\]
\end{Lemma}

\noindent {\bf Proof} 
From Lemma \ref{zetai} and the definitions of $A^i$ and $\Omega^i$ we get the following commutative diagram.

\[\begin{array}{ccccccccc}
   &             &                           &             &     0                   &             &    0
                             &             &    \\
   &             &                           &             &   \downarrow            &             &  \downarrow                                      &             &     \\
   &             &                           &             & A^i                       & \stackrel{\rho^i}{\dashrightarrow}   & \Omega^i                                           &             &    \\
   &             &                           &             &  \phantom{\beta^i} \downarrow\;_{ \beta^i}           &             &  \phantom{\mu^i} \downarrow \;_{\mu^i}
                             &             &    \\
 0 & \rightarrow & \iota^i \circ (L^i,H^{i+1}S)              & \stackrel{\nu^i}{\rightarrow} &  (H^{i+1}T,H^{i+1}S)                  & \stackrel{\eta^i}{\rightarrow} & (L^{i+1},H^{i+1}S[2])   
                             & \rightarrow & 0  \\
   &             &                           &             &   \phantom{\chi^i} \downarrow\;_{\chi^i}         &             &  \phantom{\gamma^i} \downarrow\;_{\gamma^i}                                    &             &     \\  
   &             &                           &             &  (T^{i+1},S^i[1])               & =           & (T^{i+1},S^i[1])   
                             &             &      \\
                   \end{array}        \]

The restriction  of $\eta^i$ to $A^i$ is the map $\rho^i : A^i \rightarrow \Omega^i$ such that $\rho^i \mu^i = \beta^i \eta^i.$ Hence \[Ker(\rho^i) \subseteq \iota^i \circ (L^i,H^{i+1}S).\]
On the other side,  $\iota^i \circ (L^i,H^{i+1}S) \subseteq A^i.$ In fact,  
$\nu^i \chi^i = \nu^i \eta^i \gamma^i = 0$, hence \[ \hbox{there exists} \;\;\phi^i : \iota^i \circ (L^i,H^{i+1}S) \rightarrow A^i \;\; \hbox{such that} \;\; \phi^i \beta^i = \nu^i.\] Since $\nu^i$ is injective, so is $\phi^i.$ Therefore $Ker(\rho^i) = \iota^i \circ (L^i,H^{i+1}S)$ and this finishes the proof of the lemma.\qed

By Theorem \ref{principal} and Lemma \ref{rhoi} we get  more explicit pull-backs.

\begin{Proposition} \label{expliciten} For all $i \in \{0,1,\cdots,n-2\},$ 
\[\begin{array}{ccc}
 Hom_{D^b(\Lambda)}(T^i,S^i) & \rightarrow & Hom_{D^b(\Lambda)}(T^{i+1},S^{i+1})   \\
     \downarrow           &             &  \downarrow \\
  A^i      & \stackrel{\rho^i}{\rightarrow} & \Omega^i    
                           \end{array} \] are pull-back diagrams.
\end{Proposition} 

\noindent {\bf Remark} We get morphisms \[\Sigma^i : Hom_{D^b(\Lambda)}(T,S) \rightarrow Hom_{\Lambda}(H^{i+1}T,H^{i+1}S)\] for all $i \in \{0,1,\cdots,n-2\}.$ As $Hom_{D^b(\Lambda)}(T^i,S^i) \rightarrow  Hom_{D^b(\Lambda)}(T^{i+1},S^{i+1})$ is surjective for all $i \in \{0,1,\cdots,n-2 \},$ we get $Im(\Sigma^i) = A^i.$ Hence Lemma \ref{psii} describes the cokernel of $\Sigma^i$ as a submodule of $Hom_{D^b(\Lambda)}(T^{i+1}[-1],S^i).$

\section{Torsion free $Hom$-spaces} 
 
 In this section we assume that the $R$-module $Hom_{D^b(\Lambda)}(T,S)$ is $R$-torsion free and the homologies of $T$ and $S$ are $R$-torsion modules except in degrees 1 and $n.$ We will see that $Hom_{D^b(\Lambda)}(T,S)$ is the pull-back of $A^0$ defined in section 2 and $Hom_{\Lambda}(H^nT/tors,H^nS/tors).$\\

For all $i \in \{1,\cdots,n-2\},$ we define the $R$-torsion part of $(T^i,S^i)$to be  \[t(T^i,S^i) = \{ \varphi \in (T^i,S^i) : \;\; \hbox{there exists} \;\; r \in R, r \neq 0 \;\;r.\varphi = 0 \}.\]    
 Let $\xi : (T^1,S) \rightarrow (T^1,S^1)$ be the mapping defined in Lemma \ref{thetai}.

\begin{Lemma}\label{cap} We have $t(T^1,S^1) \cap \xi((T^1,S)) = 0.$
\end{Lemma}

\noindent {\bf Proof}  Since $(T,S)$ is $R$-torsionfree, so is $(T^1,S),$ as it embeds into $(T,S).$ In particular $\xi((T^1,S))$ is also $R$-torsionfree since $\xi$ is injective, while $t(T^1,S^1)$ is $R$-torsion. \qed

We denote by $\theta : (T^1,S^1) \rightarrow \Omega^0$ the mapping defined in Proposition \ref{expliciten}, and $\theta' := \theta|_{t(T^1,S^1)}.$

\begin{Lemma}The mapping $\theta' : t(T^1,S^1) \rightarrow Im(\theta')$ is bijective. 
\end{Lemma}

\noindent {\bf Proof} The surjectivity is clear, and the injectivity results from Lemma \ref{cap}. \qed

The mapping $\theta$ induces $\hat{\theta} : \frac{(T^1,S^1)}{t(T^1,S^1)} \rightarrow \frac{\Omega^0}{Im(\theta')},$
where \[\hat{\theta}(f + t(T^1,S^1)) = \theta(f) + Im(\theta').\]

\begin{Lemma} The mapping $\hat{\theta}$ is surjective
and  $Ker(\hat{\theta}) \cong (T^1,S).$ 
\end{Lemma}

Since $\theta$ is surjective, so is $\hat{\theta}.$ It is not difficult to show that $Ker(\hat{\theta})$ is $\frac{(T^1,S) + t(T^1,S^1)}{t(T^1,S^1)}$ and by Lemma \ref{cap} we can identify this with $(T^1,S).$  \\

As a consequence, we get a pull-back diagram 
{\small
\[ \begin{array}{ccccccccc} 
& & & & 0 & & 0 & & \\
& & & & \downarrow & & \downarrow & & \\
& & & & (T^1,S) & = & \frac{(T^1,S)+t(T^1,T^S)}{t(T^1,S^1)} & & \\
& & & & \downarrow & & \downarrow & & \\
 0 & \rightarrow & t(T^1,S^1) & \rightarrow & (T^1,S^1) & \rightarrow & \frac{(T^1,S^1)}{t(T^1,S^1)} & \rightarrow & 0 \\
   &             & ||     &             & \downarrow &        & \phantom{\hat{\theta}} \downarrow \; _{\hat{\theta}}           &             &   \\
0  & \rightarrow & Im(\theta') & \rightarrow & \Omega^0 & \rightarrow & \frac{\Omega^0}{Im(\theta')} & \rightarrow & 0 \\
   &             &             &             & \downarrow &         &  \downarrow                &             &  \\
   &             &             &             &     0      &         &        0                   &             &  \end{array}\] }
\begin{equation} \label{diag1} \end{equation}
From Lemma \ref{expliciten} we have the following diagram : 
\begin{eqnarray} \label{diag2}
                        (T,S) & \rightarrow & (T^1,S^1)\nonumber \\
                        \downarrow & & \downarrow \nonumber \\
                        A^0 & \rightarrow & \Omega^0 \end{eqnarray} 
We compose the diagrams (\ref{diag1}) and (\ref{diag2}) to obtain the following diagram where the two small rectangles are pull-backs and where the kernels of the vertical mappings are isomorphic :
\[ \begin{array}{ccccc} 
                        (T,S) & \rightarrow & (T^1,S^1) & \rightarrow & \frac{(T^1,S^1)}{t(T^1,S^1)} \\
                        \downarrow & & \downarrow & & \downarrow \\
                        A^0 & \rightarrow & \Omega^0 & \rightarrow & \frac{\Omega^0}{Im(\theta')} \end{array} \]

\begin{Corollary} \label{Cor}
The composition of diagram (\ref{diag1}) and diagram (\ref{diag2})  is also a pull-back :
\[ \begin{array}{ccc} (T,S) & \rightarrow & \frac{(T^1,S^1)}{t(T^1,S^1)} \\
                      \downarrow & & \downarrow \\
                      A^0 & \rightarrow & \frac{\Omega^0}{Im(\theta')} \end{array} \]
\end{Corollary}

\begin{Lemma} \label{lem} We have $(T^1,S^1)/t(T^1,S^1) \cong (H^nT,H^nS)/t(H^nT,H^nS).$
\end{Lemma}

\noindent {\bf Proof} Fix  an $i \in \{1,\cdots,n-2\}.$ 
Let $\Phi : (T^i,S^i) \rightarrow (T^{i+1},S^{i+1})$ be the mapping defined in Lemma \ref{pii} and let $\Phi' := \Phi|_{t(T^i,S^i)}.$ 

We have the following diagram :
\[\begin{array}{ccccccccc}
  &             &                        &             & 0                    &                                 &         0                                   &             & \\
  &             &                        &             & \downarrow           &                                 &        \downarrow                           &             &  \\
  &             &                        &             & t(T^i,S^i)           & \stackrel{\Phi'}{\rightarrow} & Im(\Phi')                          &             &  \\
  &             &                        &             & \phantom{\zeta} \downarrow\; _{\zeta}          &                                 &  \phantom{\nu} \downarrow\; _{\nu}                            &             &   \\
0 & \rightarrow & \iota^i \circ(L^i,H^{i+1}S) & \stackrel{\mu}{\rightarrow} & (T^i,S^i)            & \stackrel{\Phi}{\rightarrow}  & (T^{i+1},S^{i+1})                    & \rightarrow & 0 \\
  &             &                        &             & \phantom{\lambda} \downarrow\; _{\lambda}          &                                 & \phantom{\tau}\downarrow\; _{\tau}                          &             &   \\
  &             &                        &             & (T^i,S^i)/t(T^i,S^i) & \stackrel{\delta}{\dasharrow}    &        (T^{i+1},S^{i+1})/Im(\Phi') &             & \\
  &             &                        &             & \downarrow           &                                 &         \downarrow                           &             &  \\
  &             &                        &             & 0                    &                                 &         0                                   &             &  \end{array} \]
 
  Since $\zeta \Phi \tau = \Phi' \nu \tau = 0,$ there exists \[\delta : (T^i,S^i)/t(T^i,S^i) \rightarrow (T^{i+1},S^{i+1})/Im(\Phi')\] such that
  $\Phi \tau = \lambda \delta.$ 
   Since $\Phi$ and $\tau$ are surjective, so is $\delta.$  
  
  The kernel of $\delta$ is $[\mu(\iota^i \circ(L^i,H^{i+1}S))+ t(T^i,S^i)]/t(T^i,S^i).$ Since $(H^{i+1}T,H^{i+1}S)$ is $R$-torsion, so is $\mu(\iota^i \circ(L^i,H^{i+1}S)),$ hence  $Ker(\delta) = 0.$ This implies \[(T^i,S^i)/t(T^i,S^i) \cong (T^{i+1},S^{i+1})/Im(\Phi').\] 
  
  Since $Im(\Phi') \subseteq t(T^{i+1},S^{i+1}),$ we get a surjective mapping
  \[(T^i,S^i)/t(T^i,S^i) \rightarrow (T^{i+1},S^{i+1})/t(T^{i+1},S^{i+1})\]
  whose kernel is $t(T^{i+1},S^{i+1})/Im(\Phi').$ 
  
  This kernel is formed by $R$-torsion elements while $(T^i,S^i)/t(T^i,S^i)$ is $R$-torsionfree. Hence $t(T^i,S^i) = Im(\Phi'),$ and we get
  \[(T^i,S^i)/t(T^i,S^i) \cong (T^{i+1},S^{i+1})/t(T^{i+1},S^{i+1})\]  for all $i \in \{1,\cdots,n-2\}.$
  
  As a consequence, 
  $(T^1,S^1)/t(T^1,S^1) \cong (T^{n-1},S^{n-1})/t(T^{n-1},S^{n-1}),$
  and since $T^{n-1} = Coker(\alpha^n) =  H^nT$ and $S^{n-1} = Coker(\beta^n) = H^nS,$ we obtain $(T^1,S^1)/t(T^1,S^1) \cong (H^nT,H^nS)/t(H^nT,H^nS).$ \qed

\begin{Lemma} \label{lem1} $(H^nT,H^nS)/tors \cong (H^nT/tors,H^nS/tors).$ \end{Lemma}

\noindent {\bf Proof} We denote by $L := H^nT,$ $M := H^nS,$ and by $tL$ and $tM$ the torsion parts of $L$ and $M$ respectively. We have :
\begin{enumerate} \item  $0 \rightarrow (L,tM) \rightarrow (L,M) \rightarrow (L,M/tM) \rightarrow 0$ is an exact sequence. \\
 We apply $(L,-)$ to the exact sequence $0 \rightarrow tM \rightarrow M \rightarrow M/tM \rightarrow 0$ to get 
\begin{equation} \label{LtL}
 0 \rightarrow (L,tM) \rightarrow (L,M) \rightarrow (L,M/tM) \rightarrow Ext^1(L,tM) \rightarrow \cdots
\end{equation}
Since $(Im(\alpha^n),tM) \subseteq (Im(\alpha^n),M) = 0,$ we have $Ext^1(L,tM) = 0$ and the statement holds.

 \item $(L,M/tM) \cong (L/tL,M/tM).$ \\
We apply $(-,M/tM)$ to the exact sequence $0 \rightarrow tL \rightarrow L \rightarrow L/tL \rightarrow 0$ to get
 \[ 0 \rightarrow (L/tL,M/tM) \rightarrow (L,M/tM) \rightarrow (tL, M/tM) \rightarrow \cdots \]
Since $(tL, M/tM) = 0$, we get $(L,M/tM) = (L/tL,M/tM)$. 

\item $t(L,M) = (L,tM)$ \\ 
First we need that $t(L,M) \subseteq (L,tM).$
In fact, if $\varphi \in t(L,M)$, there exists an $r \in R-\{0\},$ such that $r \varphi = 0.$ Hence there exists an $r$ such that for all $l \in L,$  $r \varphi(l) = 0.$ This implies that for all $l \in L,$  $\varphi(l) \in tM.$ 

On the other hand, if $f \in (L,tM)$ and $x \in L$,  $f(x) \in tM$ i.e there exists an $r \in R - \{0\}$ such that $rf(x) = 0.$ Since $L$ is finitely generated, there exists an $s \in R - \{0\}$ such that $sf = 0.$ Therefore, $f \in t(L,M).$ \end{enumerate}
By 1, 2, and 3, the rows of the following diagram are exact :
\[\begin{array}{ccccccccc} 0 & \rightarrow & (L,tM) & \rightarrow & (L,M) & \rightarrow & (L/tL,M/tM) & \rightarrow & 0 \\
                             &             &  ||    &             &   ||  &             &             &             &   \\
                           0 & \rightarrow & t(L,M) & \rightarrow & (L,M) & \rightarrow & (L,M)/t(L,M) & \rightarrow & 0 \end{array} \]
and this proves the lemma.
\qed

From Corollary \ref{Cor}, Lemma \ref{lem} and Lemma \ref{lem1} we get the following pull-back
  \[ \begin{array}{ccc} Hom_{D^b(\Lambda)}(T,S) & \rightarrow & Hom_{\Lambda}(H^nT/t(H^nT),H^nS/t(H^nS) \\
                      \downarrow & & \downarrow \\
                      A^0 & \rightarrow & \frac{\Omega^0}{Im(\theta')} \end{array}\]
We have just proved the following theorem

\begin{Theorem} \label{symetriquen} If
for all $i \in \{1, \cdots, n-1\},$ \[ Ext^1(L^{i+1},H^{i+1}T) = Hom_{\Lambda}(L^i,Q^i) = Hom_{\Lambda}(L^i,Im(\beta^i))= \] \[ = Hom_{\Lambda}(Im(\alpha^i),Coker(\beta^i)) = 0,\] the $R$-module $Hom_{D^b(\Lambda)}(T,S)$ is $R$-torsion free and
 $H^iT$ and $H^iS$ are $R$-torsion for all $i \in \{2,\cdots,n-2\},$ then, we get the following pull-back diagram
\[ \begin{array}{ccc} Hom_{D^b(\Lambda)}(T,S) & \rightarrow & Hom_{\Lambda}(H^nT/t(H^nT),H^nS/t(H^nS)) \\
                      \downarrow & & \downarrow \\
                      A^0 & \rightarrow & \frac{\Omega^0}{Im(\theta')} \end{array}\]
 where $A^0$ is formed by those morphisms of $Hom_{\Lambda}(H^1T,H^1S)$ that induce morphisms in $Hom_{D(\Lambda)}(T,S)$
 and $\Omega^0 = A^0/(\iota^0 \circ (P^1,H^1S)).$ \end{Theorem} 

{\bf Remarks}\begin{enumerate} \item If $R$ is a Dedekind domain, $\Lambda$ is a symmetric order, i.e. an $R$-order which is also a symmetric algebra, and  $S = T$ is a tilting complex,  then $Hom_{D^b(\Lambda)}(T,T)$ is an $R$-order by \cite[Theorem 1]{Zimmermann}, hence $R$-torsionfree.
  \item If for all $i \in \{1,\cdots,n-2\},$ the $R$-modules $Hom_{D^b(\Lambda)}(T^i,S^i)$ are $R$-torsionfree, we can get $Hom_{D^b(\Lambda)}(T^i,S^i)$ as pull-back
 \[ \begin{array}{ccc} (T^i,S^i) & \rightarrow & (T^{i+1},S^{i+1})/t(T^{i+1},S^{i+1})   \\
     \downarrow           &             &  \downarrow   \\
  A^i/tA^i      & \rightarrow & \bar{\Omega}^i \end{array} \]   
 where $A^i/tA^i$ is the $R$-torsionfree part of $A^i.$ But even though $\Lambda$ is symmetric, $Hom_{D^b(\Lambda)}(T^i,S^i)$ can be $R$-torsion, as is showed in Example \ref{exe3} with $n=2,$ $S = T$ and $T^1 = Coker(\alpha).$ 
 \item For $i = 0,$ $Ext^1(L^{i+1},H^{i+1}S) = 0$  is automatic, because we have the inclusions
$Ext^1(L^1,H^1S) \hookrightarrow (T,H^1S) \hookrightarrow (T,S)$ by Lemma \ref{Ext1i} and Lemma \ref{psii}, and $(T,S)$ is an $R$-torsionfree $\Lambda$-module. \end{enumerate}              

\section{Examples}

We will give some examples demonstrating how to use Theorem \ref{principal}.  We use rings $\Lambda$ which are orders and we take $S = T.$  The presentations of $\Lambda$ in Example \ref{exe3} and Example \ref{exe4} can be found in \cite{Roggenkamp} or in \cite[section 4.4]{Koenig-Zimmermann2}. 

We abbreviate 
$(X,Y) := Hom_{D^b(\Lambda)}(X,Y)$ and $Ext^1(X,Y) := Ext^1_{\Lambda}(X,Y)$ for two complexes $X$ and $Y.$

Let $K$ be $Frac(R)$ and $\Lambda$ be an $R$-order. Let $\chi_1, \cdots, \chi_m$ be the irreducible characters of $K \otimes_R \Lambda.$ Then we write $L_i$ for a $\Lambda$-lattice $L$ such that $K \otimes_R \Lambda$ acts on $K \otimes_R L$ as $\chi_i.$
 
\begin{Exe} \label{exe3} Let $R$ be the 5-adic integers  $\mathbb{Z}_5,$ and $\pi = Rad(R) = <5>.$ Let $\Lambda$ be the $R$-order
\[\left \{ \left( d_0, \left( \begin{array}{cc} a_1 & b_1 \\ c_1 & d_1 \end{array} \right),   \left( \begin{array}{cc} a_2 & b_2 \\ c_2 & d_2 \end{array} \right),   \left( \begin{array}{cc} a_3 & b_3 \\ c_3 & d_3 \end{array} \right), a_4 \right), a_i,b_i,c_i,d_i \in R, 5|(d_i-a_{i+1}), 5|c_i \right\} \] which we write (see \cite{Roggenkamp} or \cite[section 4.4]{Koenig-Zimmermann2}) as   

\begin{center}
\unitlength1cm
\begin{picture}(8,1)
\put(1,0){$R_1$}
\put(1.25,0.25){\line(3,1){0.75}}
\put(1.5,0){$\left( \begin{array}{cc} R & R \\
                  \pi & R   \end{array} \right)_2$}
\put(3,0){\line(3,1){1}}
\put(3.5,0){$\left( \begin{array}{cc} R & R \\
                  \pi & R   \end{array} \right)_3$} 
\put(5,0){\line(3,1){1}}                 
\put(5.5,0){$\left( \begin{array}{cc} R & R \\
                  \pi & R   \end{array} \right)_4$} 
\put(7,0){\line(3,1){0.5}}                   
\put(7.5,0){$R_5,$}                                                    
\end{picture}
\end{center}  
\phantom{l}
where the indices indicate the characters involved in the indecomposable projective modules. Then, by \cite{Roggenkamp}, $\Lambda$ is Morita equivalent to $B_0(\mathbb{Z}_5\mathcal{S}_5),$ the block principal of the group ring of the symmetric group of degree 5 over the 5-adic integers.  \\
Let
\begin{center}
\unitlength1cm
\begin{picture}(8,0)
\put(0,0){$P_0 =$}
\put(1,0){$R_1$}
\put(1.25,0.25){\line(3,1){0.75}}
\put(1.5,0){$\left( \begin{array}{c} R  \\
                  \pi   \end{array} \right)_2,$}
\put(3.5,0){$P_1 = $}
\put(4.5,0){$\left( \begin{array}{c} R  \\
                  R   \end{array} \right)_2$}
\put(5.5,0){\line(3,1){1}}
\put(6,0){$\left( \begin{array}{c} R  \\
                  \pi   \end{array} \right)_3,$}
\end{picture}
\end{center}
\phantom{l}
\begin{center}
\unitlength1cm
\begin{picture}(8,0)
\put(0,0){$P_2 =$}
\put(1,0){$\left( \begin{array}{c} R  \\
                  R   \end{array} \right)_3$}
\put(2,0){\line(3,1){1}}
\put(2.5,0){$\left( \begin{array}{c} R  \\
                  \pi   \end{array} \right)_4,$}
                  \put(4.5,0){$P_3 = $}
\put(5.5,0){$\left( \begin{array}{c} R  \\
                  R   \end{array} \right)_4$}
\put(6.5,0){\line(3,1){0.5}}
\put(7,0){$R_5.$}
\end{picture}
\end{center}
By definition, the characters of $P_0$ are $\chi_1$ and $\chi_2$, those of $P_1$ are $\chi_2$ and $\chi_3$, those of $P_2$ are $\chi_3$ and $\chi_4,$ and those of $P_3$ are $\chi_4$ and $\chi_5.$

Let $T$ be the complex $T_1 \oplus T_2 \oplus T_3 :$
 \[ \begin{array}{ccccccccc} T_1 &:& 0 & \rightarrow & P_1 \oplus P_3 & \stackrel{\alpha}{\rightarrow} & P_2&  \rightarrow & 0 \\
                           \oplus     & &   &             & \oplus &                             &   &              &   \\                                T_2 &:& 0 & \rightarrow & P_1 \oplus P_3 & \rightarrow                    & 0  &              & \\               \oplus  & &   &             & \oplus &                             &    &              &   \\  
                            T_3 &:& 0 & \rightarrow & P_0            & \rightarrow & 0 & & \end{array} \]
The complex $T$ is a tilting complex according to  the proof of \cite[Lemma 5.1.2]{Koenig-Zimmermann2}. 

We define $\bar{\alpha} = (\alpha,0,0)$ and observe that $End_{\Lambda}(Ker(\bar{\alpha}))$ decomposes as a direct product of two rings :
\begin{center}
\unitlength1cm
\begin{picture}(8,2)
\put(0.75,0){$R_4$}
\put(1.25,0.25){\line(3,1){0.75}}
\put(1.5,0){$\left( \begin{array}{cc} R & R \\
                  \pi & R   \end{array} \right)_5$}
\put(3.5,0){$\oplus$}
\put(4,0){$\left( \begin{array}{ccc} R & R & \pi \\
                  \pi & R  & \pi \\ 
                  R & R & R \end{array} \right)_2$} 
\put(6,1.5){$R_3$}
\put(6,1.5){\line(-1,-2){0.6}}
\put(7,0){\line(-3,-1){0.75}}                   
\put(7,0){$R_1$}                                                    
\end{picture}
\end{center} 
\phantom{k}
 If the mapping $\psi$ defined in Lemma \ref{psii} was surjective,then since $End_{\Lambda}(Coker(\bar{\alpha}))$ is an $R$-torsion module, we would get $End_{D^b(\Lambda)}(T) \cong End_{\Lambda}(Ker(\bar{\alpha}))$ which is decomposable. This is impossible. Hence, $\psi$ cannot be surjective.  

In fact, $End_{D^b(\Lambda)}(T) \cong A^0$ is the ring :
\begin{center}
\unitlength1cm
\begin{picture}(8,2)
\put(0.75,0){$R_4$}
\put(1.25,0.25){\line(3,1){0.75}}
\put(1.5,0){$\left( \begin{array}{cc} R & R \\
                  \pi & R   \end{array} \right)_5$}
\put(3,0){\line(3,1){1.25}}
\put(4,0){$\left( \begin{array}{ccc} R & R & \pi \\
                  \pi & R  & \pi \\ 
                  R & R & R \end{array} \right)_2$} 
\put(6,1.5){$R_3$}
\put(6,1.5){\line(-1,-2){0.6}}
\put(7,0){\line(-3,-1){0.75}}                   
\put(7,0){$R_1$}                                                    
\end{picture}
\end{center} 
\phantom{t} 
\phantom{u}  
where the congruence $R_2 - R_5$ is a consequence of the congruence in $P_2.$ 
\end{Exe} 

\begin{Exe} \label{exe4}Let $R$ be the 7-adic integers $\mathbb{Z}_7,$ the 7-adic integers,  $\pi = Rad(R) = <7>,$ and let $\Lambda$ be the  $R$-order 
\begin{center}
\unitlength1cm
\begin{picture}(12,1)
\put(1,0){$R$}
\put(1.25,0.25){\line(3,1){0.75}}
\put(1.5,0){$\left( \begin{array}{cc} R & R \\
                  \pi & R   \end{array} \right)$}
\put(3,0){\line(3,1){1}}
\put(3.5,0){$\left( \begin{array}{cc} R & R \\
                  \pi & R   \end{array} \right)$} 
\put(5,0){\line(3,1){1}}                 
\put(5.5,0){$\left( \begin{array}{cc} R & R \\
                  \pi & R   \end{array} \right)$} 
\put(7,0){\line(3,1){1}}                
\put(7.5,0){$\left( \begin{array}{cc} R & R \\
                  \pi & R   \end{array} \right)$}  
\put(9,0){\line(3,1){1}}                
\put(9.5,0){$\left( \begin{array}{cc} R & R \\
                  \pi & R   \end{array} \right)$}                                       
\put(11,0){\line(3,1){0.5}}                   
\put(11.5,0){$R$}                                                    
\end{picture}
\end{center}  
\phantom{l}
By \cite{Roggenkamp} $\Lambda$ is Morita equivalent to $B_0(\mathbb{Z}_7\mathcal{S}_7),$ the block principal of the group ring of the symmetric group of degree 7 over the 7-adic integers. \\

Let $T$ be the complex 
\[ \begin{array}{ccccccccc} 
 0 & \rightarrow & P_0  \oplus P_0          & \stackrel{(\alpha,0)}{\rightarrow} & P_1             &                   \rightarrow                   & 0   &             & \\
   &             & \oplus         &                                & \oplus          &                                                      & \oplus    &             & \\
   &             & 0              & \rightarrow                    & P_2 \oplus P_4  & \stackrel{\beta}{\rightarrow} & P_3    & \rightarrow & 0 \\
   &             & \oplus         &                                &   \oplus              &                                                      & \oplus &             & \\
   &             & 0              & \rightarrow                    & P_2 \oplus P_4 \oplus P_5 &                        \rightarrow                   & 0      &             &  \end{array} \]              
where $P_1,$ $P_2$, $P_3,$ $P_4$ and $P_5$ are projective $\Lambda$-modules :
\begin{center}
\unitlength1cm
\begin{picture}(8,0)
\put(0,0){$P_0 =$}
\put(1,0){$R_1$}
\put(1.25,0.25){\line(3,1){0.75}}
\put(1.5,0){$\left( \begin{array}{c} R  \\
                  \pi   \end{array} \right)_2,$}
\put(3.5,0){$P_1 = $}
\put(4.5,0){$\left( \begin{array}{c} R  \\
                  R   \end{array} \right)_2$}
\put(5.5,0){\line(3,1){1}}
\put(6,0){$\left( \begin{array}{c} R  \\
                  \pi   \end{array} \right)_3,$}
\end{picture}
\end{center}
\phantom{r}
\begin{center}
\unitlength1cm
\begin{picture}(8,0)
\put(0,0){$P_2 =$}
\put(1,0){$\left( \begin{array}{c} R  \\
                  R   \end{array} \right)_3$}
\put(2,0){\line(3,1){1}}
\put(2.5,0){$\left( \begin{array}{c} R  \\
                  \pi   \end{array} \right)_4,$}
\put(4.5,0){$P_3 = $}
\put(5.5,0){$\left( \begin{array}{c} R  \\
                  R   \end{array} \right)_4$}
\put(6.5,0){\line(3,1){1}}
\put(7,0){$\left( \begin{array}{c} R  \\
                  \pi   \end{array} \right)_5,$}
\end{picture}
\end{center}
\phantom{k}
\begin{center}
\unitlength1cm
\begin{picture}(8,0)
\put(0,0){$P_4 =$}
\put(1,0){$\left( \begin{array}{c} R  \\
                  R   \end{array} \right)_5$}
\put(2,0){\line(3,1){1}}
\put(2.5,0){$\left( \begin{array}{c} R  \\
                  \pi   \end{array} \right)_6,$}
\put(4.5,0){$P_5 = $}
\put(5.5,0){$\left( \begin{array}{c} R  \\
                  R   \end{array} \right)_6$}
\put(6.5,0){\line(3,1){0.5}}
\put(7,0){$R_7,$}
\end{picture}
\end{center}
\phantom{l}
and the indices indicate the characters involved in the modules. 

We define by $\bar{\alpha} := ((\alpha,0),0,0)$ and $\bar{\beta} := (0,\beta,0).$ \\
Since $(Im(\bar{\alpha}),Ker(\bar{\alpha})) = (Im(\bar{\beta}),H^2T) = 0,$ we get \[Ext^1(Coker(\bar{\alpha}),Ker(\bar{\alpha})) = Ext^1(Coker(\bar{\beta}),H^2T) = 0.\]

Let $T^1$ be the complex \[ 0 \rightarrow Im(\bar{\alpha}) \rightarrow \left( \begin{array}{c} P_1 \\ P_2 \oplus P_4 \\ P_2 \oplus P_4 \oplus P_5  \end{array} \right) \rightarrow \left( \begin{array}{c} 0  \\ P_3  \\ 0 \end{array} \right) \rightarrow 0.\]
By Theorem \ref{principal}, $End_{D^b(\Lambda)}(T^1)$ is a pull-back of $End_{\Lambda}(H^2T)$ and  \\$End_{\Lambda}(Coker(\bar{\beta}))$ over $Ext^2(Coker(\bar{\beta}),H^2T).$ \\

The mapping $(Coker(\bar{\beta}),Coker(\bar{\beta})) \stackrel{\theta^1}{\rightarrow} Ext^2(Coker(\bar{\beta}),H^2T)$ is induced by $Ker(\bar{\beta}) \rightarrow H^2T$ and this is  multiplication by an $r \in R.$ 

Hence, $\theta^1$ is injective, and \[End_{D^b(\Lambda)}(T^1) \cong Im(\psi^1)\] where $\psi^1 : (T^1,T^1) \rightarrow (H^2T,H^2T).$  In fact, $A^1 = Im(\psi^1)$ is the following ring :
\begin{center}
\unitlength1cm
\begin{picture}(10,2)
\put(1.75,0){$\left( \begin{array}{ccc} R & R & R \\
                  \pi & R  & R \\ 
                  \pi & \pi & R \end{array} \right)_3$} 
\put(3.8,1.5){$R_4$}
\put(3.8,1.5){\line(-1,-2){0.6}}
\put(5.2,0.5){\line(-3,-2){1.25}}
\put(4.75,0){$\left( \begin{array}{ccc} R & R & R \\
                  \pi & R  & R \\ 
                  \pi & \pi & R \end{array} \right)_6$}
\put(6.75,1.5){$R_5$}
\put(6.75,1.5){\line(-1,-2){0.6}}                  
\put(7.6,0){\line(-3,-1){0.75}}                   
\put(7.75,0){$R_7$}                                                    
\end{picture}
\end{center} 
\phantom{l}
where the congruence of $R_3 - R_6$ comes from the one in $P_3.$ \\
 
 Now, $End_{D^b(\Lambda)}(T)$ is a pull-back of $End_{\Lambda}(Ker(\bar{\beta}))$ and $End_{D^b(\Lambda)}(T^1)$ over $Ext^2(Coker(\bar{\alpha}),Ker(\bar{\alpha})).$
 
The only non-zero component of $Ext^2(Coker(\bar{\alpha}),Ker(\bar{\alpha}))$ is $R/\pi,$ which gives the congruence $R_1 - R_3.$ Hence, $End_{D^b(\Lambda)}(T)$ is the following ring : 
\begin{center}
\unitlength1cm
\begin{picture}(10,2)
\put(0.5,0){$R_2$}
\put(1,0){\line(3,1){1}}
\put(1.5,0){$\left( \begin{array}{cc} R & R \\
                  \pi & R   \end{array} \right)_1$}
\put(3,0){\line(3,1){1.25}}
\put(3.75,0){$\left( \begin{array}{ccc} R & R & R \\
                  \pi & R  & R \\ 
                  \pi & \pi & R \end{array} \right)_3$} 
\put(5.8,1.5){$R_4$}
\put(5.8,1.5){\line(-1,-2){0.6}}
\put(7.2,0.5){\line(-3,-2){1.25}}
\put(6.75,0){$\left( \begin{array}{ccc} R & R & R \\
                  \pi & R  & R \\ 
                  \pi & \pi & R \end{array} \right)_6$}
\put(8.75,1.5){$R_5$}
\put(8.75,1.5){\line(-1,-2){0.6}}                  
\put(9.6,0){\line(-3,-1){0.75}}                   
\put(9.75,0){$R_7$}                                                    
\end{picture}
\end{center} 
\phantom{f}

{\bf Remarks} \begin{enumerate} \item The complex $T$ is a tilting complex. If $\psi^1 :(T^1,T^1) \rightarrow (H^2T,H^2T)$ was  surjective, we would get  $(T^1,T^1) \cong (H^2T,H^2T)$ which is decomposable and this would impliy that $(T,T)$ is decomposable, which is a contradiction. Hence, in this case, $\psi^1$ cannot be surjective.
\item Nevertheless, the mapping $\psi^0 : (T,T) \rightarrow (H^1T, H^1T)$ is surjective.
\end{enumerate}
\end{Exe}

{\footnotesize

\noindent Address : \\
{\sc Facult\'e de Math\'ematiques et CNRS (LAMFA UMR 6140), Universit\'e de Picardie, 33 rue St Leu, 80039 Amiens Cedex, France} \\
 {\it Email address :} intan.muchtadi@u-picardie.fr \\

\noindent Current address : \\
{\sc Department of Mathematical Sciences, NTNU, N-7491 Trondheim Norway} \\
{\it Email address :} alamsyah@math.ntnu.no \\

\noindent Permanent address : \\
{\sc Departemen Matematika, Institut Teknologi Bandung (ITB), Jl. Ganesha no. 10, Bandung 40132, Indonesia} \\
{\it Email address :} ntan@dns.math.itb.ac.id}

\end{document}